\documentclass[conference]{IEEEtran}

\ifCLASSINFOpdf
\else
\fi

\usepackage[noadjust]{cite}
\usepackage[caption=false,font=footnotesize]{subfig}
\usepackage{amsmath,amssymb,amsfonts}
\usepackage{algorithmic}
\usepackage{graphicx}
\usepackage{bbm}
\usepackage{textcomp}
\usepackage{xcolor}
\usepackage{url}

\usepackage[]{fancyhdr} %
\newcommand{\changefont}{\fontsize{7}{7}\selectfont}
\newcommand{\up}[1]{^\mathrm{#1}}

\IEEEoverridecommandlockouts
\begin{document}

\title{Impact of Bidding and Dispatch Models over Energy Storage Utilization in Bulk Power Systems}

\author{\IEEEauthorblockN{Ningkun Zheng,  Bolun Xu}
\IEEEauthorblockA{Earth and Environmental Engineering \\Columbia University\\
New York, NY, USA\\
\{nz2343, bx2177\}@columbia.edu}
}

\maketitle
\thispagestyle{fancy}
\pagestyle{fancy}


\begin{abstract}

Energy storage is a key enabler towards a low-emission electricity system, but requires appropriate dispatch models to be economically coordinated with other generation resources in bulk power systems.
This paper analyzes how different dispatch models and bidding strategies would affect the utilization of storage with various durations in deregulated power systems. We use a dynamic programming model to calculate the operation opportunity value of storage from price predictions, and use the opportunity value result as a base for designing market bids. We compare two market bidding and dispatch models in single-period economic dispatch: one without state of charge (SoC) constraints and one with SoC constraints. We test the two storage dispatch models, combined with different price predictions and storage durations, using historical real-time price data from New York Independent System Operator. We compare the utilization rate with respect to results from perfect price forecast cases. Our result shows that while price prediction accuracy is critical for short duration storage with a less than four hours capacity, storage with a duration longer than twelve hours can easily achieve a utilization rate higher than 80\% even with naive day-ahead price predictions. Modeling storage bids as dependent of SoC in single-period real-time dispatch will provide around 5\% of improvement in storage utilization over all duration cases and bidding strategies, and higher renewable share will likely improve storage utilization rate due to higher occurrence of negative prices.

\end{abstract}



%
\IEEEpeerreviewmaketitle

\section{Introduction}
A decarbonizing electricity system based on renewable energy resources  faces challenges in system reliability and economical efficiency. With the ambitious renewable portfolio standard and carbon emission limitation targets,  renewable capacities could surge exponentially to meet electricity demand and regional decarbonization targets~\cite{shaner2018geophysical,sepulveda2018role}. Energy storage with a less than four hours duration help to flatten the daily demand curve and renewable fluctuations, reduce the cost of electricity and mitigate the occurrence of extreme price spikes~\cite{lee2022targeted}. Long-duration energy storage has the potential to store large quantities of renewable energy over a long time scale, carrying low-cost energy from a high renewable availability period to a renewable "drought" season  is valuable to lower the system cost and reduce carbon emission. The U.S. Department of Energy aims at developing economical utility-scale long-duration energy storage within one decade, reducing the production cost by 90\% compared to current lithium-ion batteries for energy storage with more than 10 hours of duration~\cite{doe}.

Currently, lithium-ion batteries and flow batteries are the two dominant technology groups for building energy storage with one to six hours of duration. Longer duration battery energy storage is technically feasible, but economically undesirable due to the relatively high cost of storage capacity~\cite{zakeri2015electrical}. Pumped hydro storage (PHS) and compressed air energy storage (CAES) take on the function of sustaining electricity supply during multi-day periods of average demand exceeding average renewable supply~\cite{safaei2015much}. Although these mechanical energy storage technologies are relatively mature, they suffer from several limitations and may not be the optimal choice for long-duration energy storage~\cite{jenkins2021long}. The PHS and CAES are both geographically constrained from sitting in the desired location~\cite{sepulveda2021design}, while the PHS has more environmental concerns and the large-scale CAES designs combust non-renewable natural gas ~\cite{kucukali2014finding,blanco2018review,fertig2011economics}. There is a multitude of possible long-duration energy storage technologies available, including electrochemical, chemical, thermal, and previously mentioned mechanical options~\cite{sepulveda2021design}. Each technology has its own set of cost and performance characteristics. 

While the future technological pathway of energy storage is under rapid development,  integration of storage in the power grids is of guiding significance to efficiently integrate renewable and reduce electricity costs.
Power system operators must integrate energy storage into power system dispatch to economically coordinate storage  with other grid resources. Unlike conventional thermal generators whose operating costs solely depend on fuel prices and their heat rate efficiency curves,  the operation of energy storage is critically dependent on future opportunities for discharging and charging the storage. Energy storage opportunity value is highly correlated to volatilities in the system, which are jointed contributed by the share of renewables~\cite{guerra2020value}, daily demand variations with slight weekly patterns~\cite{conejo2005day}, and operating constraints from thermal generators including ramp rates, start-up and shut-down limits~\cite{kirschen2018fundamentals}. Many of these variation factors require operational planning that spans days or weeks to capture their impact over storage operation opportunities, such that storage could effectively shift energy through time to improve renewable integration and reduce total system operating cost~\cite{riddervold2021internal}. In deregulated power systems, optimized operation planning will also provide improved revenue to storage participants by dispatching storage to charge during price valleys and discharge during price spikes~\cite{li2014economic, o2017efficient}.

Existing power system dispatch models consider limited future information in the optimization. System operators must finish the dispatch optimization problem and clear the market within a fixed time frame, around 20 minutes in day-ahead markets, and a few minutes in real-time markets~\cite{chen2022battery}. Computation complexity is often the limiting factor for system operators to further increase the planning horizon. Day-ahead unit commitments usually consider a 24-hour horizon and optimize storage operation over daily patterns~\cite{huang2021configuration}. However, further extending the horizon for integrating long-duration energy storage may cause computation challenges and require significant changes to the market clearing procedure. Real-time economic dispatches usually consider a single time period or with a limited look-ahead of around one hour for ramping limits and real-time unit commitments~\cite{nyiso_rtd,caiso_rtd}. While studies have explored real-time market design using multi-period economic dispatch~\cite{guo2021pricing}, extending the real-time dispatch optimization horizon to 24 hours or more is unlikely to be feasible in practice due to computation complexities~\cite{zhao2019multi}. An alternative option is to help storage self-schedule~\cite{conejo2002self}, but this will eliminate the possibility to utilize storage in mitigating uncertainties in real-time power system operations such as renewable fluctuations or demand surges.

Facilitated by FERC (Federal Energy Regulatory Commission) Order 841~\cite{ferc_order_841}, all system operators in the United States have implemented  storage bidding models in wholesale electricity markets, which allow storage to bid as a combination of generator and flexible demand~\cite{sakti2018review}. In this model, storage participants submit a discharge bid indicating a willingness to discharge if the market clearing price is higher than the bid, and a charge bid indicating a willingness to charge if the market clearing price is lower than the bid. Therefore, storage participants engineer their discharge and charge opportunities into bid values. System operators can dispatch storage without explicitly modeling their state-of-charge limits and future system conditions, which only require marginal modification to current dispatch software. Yet, stakeholders have expressed a desire to more effectively manage their storage resource’s SoC in the real-time dispatch~\cite{caiso_es}, but explicit modeling SoC constraints require significant modification of the dispatch software to include sufficiently long predictive look-ahead horizon, which also comes with computation challenges as the joint effect of extending optimization horizon and SoC constraints may exponentially increase the computation time~\cite{chen2022battery}.

While many studies have investigated strategic storage participation in wholesale electricity markets with the objective to maximize market revenue with the potential to exercise market power~\cite{wang2017look, shafiee2016risk, krishnamurthy2017energy,thatte2013risk, salles2017potential}. The focus of this paper is to investigate how different market bidding and dispatch models would impact the utilization of storage in terms of market revenue and efficiency to reduce total system costs. Motivated to minimize the computation complexity and software upgrades needed to incorporate storage in existing bulk power system dispatch models, we only consider dispatch models without inter-temporal constraints. These models can be directly incorporated into real-time power system dispatches, and into day-ahead unit commitments to reduce the computation complexity that arose from modeling SoC constraints. This paper proposes a bidding design method using the opportunity value of energy storage, and use it explores the utilization of energy storage of various durations in bulk power system dispatch with different dispatch models.  The contribution of this paper is listed as follows
\begin{enumerate}
    \item We use a dynamic programming model to accurately calculate the opportunity value of energy storage using predicted future price data and physical characteristics of the storage including discharge cost, efficiency, and energy storage duration.
    \item We consider two types of bidding model for energy storage in single-period power system dispatch: a power bid model in which storage submits bids for charge and discharge, and an SoC bid model in which storage submits piece-wise linear bids dependent on its SoC level. Both models consider a single-period dispatch setting and do not need inter-temporal constraints to optimize storage operation.
    \item We propose a bidding strategy to design the power bids and SoC bids using the storage valuation results from the dynamic programming valuation method. 
    \item We simulate storage operation using the proposed dispatch and bidding model with historical price data from New York Independent System Operator. We compare the utilization of storage in terms of market revenue from different storage durations and dispatch models. 
\end{enumerate}

The rest of the paper is organized as follows: Section~II introduces models for storage valuation, dispatch, and bidding. Section~III introduces computation experiment settings and data sources. Section~IV presents results and Section~V concludes the paper.

\section{Models and Methods}
We introduce two types of single-period real-time dispatch models for energy storage: a power bid model and an SoC bid model. In the power bid model, the storage submits charge and discharge bids, and the dispatch decision is independent of the storage SoC.  In the SoC bid model, the bids are dependent on the storage SoC; hence the system operator will update the SoC of all storage participants at the beginning of each dispatch period, and model the discharge and charge cost depending on storage SoC.

We introduce an analytical method to calculate the opportunity value of energy storage using a price arbitrage formulation, and use the result to design bids into the two considered storage dispatch models. The opportunity value is based on predictions of future electricity prices and the physical parameter of the storage including discharge cost, efficiency, and duration. We first present a formulation of dynamic programming for energy storage arbitrage, followed by an analytical method to calculate opportunity value as a function of the SoC level. We then introduce the two dispatch models and explain how to engineer the opportunity value function into market bids.

\subsection{Storage Valuation using Dynamic Programming}
We formulate the energy storage opportunity valuation as a price response arbitrage problem using dynamic programming. We assume energy storage is a price taker and the  arbitrage problem is formulated as
\begin{subequations}\label{eq1}
\begin{align}
    Q_{t-1}(e_{t-1}) &= \max_{b_t, p_t} \lambda_t (p_t-b_t) - cp_t + Q_{t}(e_{t}) 
    \label{eq:obj2}
\end{align}
where $Q_{t-1}$ is the maximized energy storage arbitrage profit dependent on the energy storage SoC at the end of the previous time period $e_{t-1}$. This profit accounts from time period $t$ till the end of the optimizing horizon $T$. The energy market revenue is the product of the real-time market price $\lambda_t$ and the energy storage dispatch decision $(p_t-b_t)$, where $p_t$ is the discharge power and $b_t$ is the charge power. The discharge cost is the second term of the objective function, where $c$ is the marginal discharge cost. $Q_t$ represents the opportunity value of the energy storage SoC $e_t$ at the end of time period $t$, hence the value-to-go function in dynamic programming. 

The objective function subjects to the following constraints
\begin{gather}
    0 \leq b_t \leq P,\; 0\leq p_t \leq P \label{p1_c2} \\
    \text{$p_t = 0$ if $\lambda_t < 0$} \label{p1_c5}\\
    e_t - e_{t-1} = -p_t/\eta + b_t\eta \label{p1_c1}\\
    0 \leq e_t \leq E \label{p1_c3}
\end{gather}
\end{subequations}
where \eqref{p1_c2} models the upper bound $P$ and lower bound 0 of storage charge and discharge power. \eqref{p1_c5} is a relaxed form of constraint that enforces the energy storage to not discharge and charge simultaneously. Negative price is the necessary condition for storage to charge and discharge at the same time in price arbitrage, hence by enforcing storage not to discharge when the price is negative, we eliminate simultaneous charge and discharge~\cite{xu2019operational}. \eqref{p1_c1} models the energy storage SoC evolution constraint with efficiency $\eta$. \eqref{p1_c3} models the upper bound $E$ and lower bound (we assume as 0) of the storage SoC level.

\subsection{Solution Algorithm}
We define $q_t$ as the derivative of storage opportunity value function $Q_t$ in \eqref{eq:obj2}, which represents the marginal opportunity value of energy stored in the storage. It is evidently that $Q_t$ is valued and differentiable over the energy storage SoC level $e_t$ within [0\,,\,$E$]. We can move the derivative operation into the expectation calculation in \eqref{eq:obj2}. Then we can use an analytical formulation to calculate the opportunity value $q_t(e)$ at any given energy storage SoC level.

Our prior work~\cite{xu2019operational} proved $q_{t-1}$ can be recursively calculated with next period value function $q_t$, power rating $P$, and efficiency $\eta$. We rewrite this value function calculating using the deterministic formulation investigated in this paper as
\begin{align}\label{eq3}
    &q_{t-1}(e) = \nonumber\\
    &\begin{cases}
    q_{t}(e+P\eta)  & \text{if $\lambda_{t}\leq q_{t}(e+P\eta)\eta$} \\
    \lambda_{t}/\eta  & \text{if $ q_{t}(e+P\eta)\eta < \lambda_{t} \leq q_{t}(e)\eta$} \\
    q_{t}(e) & \text{if $ q_{t}(e)\eta < \lambda_{t} \leq [q_{t}(e)/\eta + c]^+$} \\
    (\lambda_{t}-c)\eta & \text{if $ [q_{t}(e)/\eta + c]^+ < \lambda_{t}$} \\
    & \quad\text{$ \leq [q_{t}(e-P/\eta)/\eta + c]^+$} \\
    q_{t}(e-P/\eta) & \text{if $\lambda_{t} > [q_{t}(e-P/\eta)/\eta + c]^+$} 
    \end{cases}
\end{align}
which calculates the opportunity value function assuming the price follows a recursive computation framework. Thus we are able to get opportunity value function $q_t(e)$ at any time period using backward calculation by stating an end period value function $q_T$.

We further discretize $q_t$ by equally spaced energy storage SoC level $e$ into small segments, which is far smaller than power rating $P$. For any SoC level $e_t$, we can find the nearest segment and return the corresponding value. Note that $Q_t$ in objective function is the integral of $q_t$. Therefore, discretizing the derivative $q_t$ is equivalent to approximate $Q_t$ using piece-wise linear functions.

\subsection{Economic Dispatch with Storage Power Bids}
We consider the existing storage economic dispatch model in which dispatches storage using charge and discharge bids, and engineer the storage marginal opportunity value $q_t(e)$ into bidding design. In this model, a storage participant submits two hourly bids, a discharge bid specifying the price above which the storage is willing to discharge, and a charge bid specifying the price below which the storage is willing to charge. Since the storage dispatch decision entirely depends on the submitted bids, the system operator is relieved from modeling the SoC constraint. We call this model the \emph{power bid} model as the bids are directly associated with the storage charge and discharge power.

The single-period economic dispatch model with storage power bids over dispatch period $t$, including $k\in K$ generators and $m\in M$ storage participants, is formulated as
\begin{subequations}\label{ed1}
\begin{align}
    &\min_{p_{m,t}, b_{m,t}, g_{k,t}} \sum_{k} C_k(g_{k,t}) + \sum_{m}c\up{p}_{m,t} p_{m,t} - c\up{b}_{m,t} b_{m,t} \label{ed1:obj}\\
    &\text{subjects to }\nonumber\\
    &0 \leq p_{m,t} \leq P\up{p}_m(e_{m,t-1})  \label{ed1:c1}\\
    &0 \leq b_{m,t} \leq P\up{b}_m(e_{m,t-1})  \label{ed1:c2}\\
    & \text{other generator, network, nodal balance constraints}\nonumber
\end{align}
\end{subequations}
Note that we keep the subscript $t$ to all decision variables but the problem considers only one time period, hence this problem represents a sequential single-period economic dispatch problem, in which the system operator optimizes a single period at a time, update dispatch decisions, and proceed to the next time period.
Decisions variables include the discharge $p_{m,t}$ and charge $b_{m,t}$ dispatch for storage $m$ over time period $t$, and generation dispatch $g_{k,t}$ for generators.  $C_k(\cdot)$ is the production cost function for generator $k$, $c\up{p}_m$ and $c\up{b}_{m,t}$ are the discharge and charge bids for storage $m$, respectively. Note that storage bids have subscript $t$ to denote these bids are time-dependent, i.e., storage will bid different values over different time periods. The problem subjects to storage discharge and charge power limits \eqref{ed1:c1} and \eqref{ed1:c2}, note that we represent the power rating as a function on the storage SoC $e_{t-1}$ at the beginning of the dispatch period, such that the dispatch decision will not order the storage to over charge or over discharge itself. We consider $e_{t-1}$ as a parameter given to the economic dispatch instead of a decision variable. This constraint follows recent market designs that the system operator will monitor storage SoC and update their power ratings in the dispatch software accordingly~\cite{caiso_es}. The economic dispatch includes other typical generator, network, and nodal power balance constraints, which we omit in the formulation for simplicity.

We now design discharge and charge bids using the opportunity valuation results. Note that these bids represent the combination of the discharge cost and the change in the opportunity value. Since the bids do not contain information about storage SoC, we assume equal possibilities of storage SoC at the beginning of the dispatch period, hence we calculate the discharge bid value as the average of marginal discharge cost with respect to battery SoC as 
\begin{align}
    c\up{p}_t &= \frac{1}{E} \int_{0}^E \frac{\partial}{\partial p_t} (c p_t - Q_t(e_{t-1} - p_t/\eta + b_t\eta)) d{e_{t-1}}\nonumber\\
    &= c+ \frac{1}{E} \int_{0}^E q_t(e_{t-1} - p_t/\eta + b_t\eta) d{e_{t-1}} /\eta \nonumber\\
    &\approx c+ \frac{1}{\eta E} \int_{0}^E q_t(e) d{e} \nonumber\\
    &= c + \bar{q}_t/\eta
\end{align}
where $\bar{q}_t$ is the average of the marginal opportunity value of storage over all SoC levels
\begin{align}
    \bar{q}_t = \frac{1}{E} \int_{0}^E q_t(e) d{e}
\end{align}
hence we derive the discharge bid as the sum of the storage discharge cost and the average marginal opportunity value $\bar{q}_t$ divided by the discharge efficiency. Similarly, for the charge bid
\begin{align}
    c\up{b}_t &= \frac{1}{E} \int_{0}^E \frac{\partial}{\partial b_t} (c p_t - Q_t(e_{t-1} - p_t/\eta + b_t\eta)) d{e_{t-1}}\nonumber\\
    &= \frac{1}{E} \int_{0}^E q_t(e_{t-1} - p_t/\eta + b_t\eta) d{e_{t-1}} \eta \nonumber\\
    &\approx \frac{\eta}{E} \int_{0}^E q_t(e) d{e} \nonumber\\
    &= \bar{q}_t\eta
\end{align}
which is the average marginal opportunity value multiplied by charging efficiency. Fig.~\ref{fig:soc_bid1} shows an example of discharge and charge bids and the marginal SoC opportunity value $q_t(e)$ used to calculate the bids. 

\subsection{Economic Dispatch with Storage SoC Bids}
Motivated by the piece-wise linear heat rate curve bids for conventional generators, we envision a market model in which the system operator lets storage design market bids as piece-wise linear functions to its SoC. The system operator still performs a single-period economic dispatch, but will add constraints to model the change in the storage SoC and the opportunity value over the dispatch period. Different from the power bids model, in this case the storage participants submit bids segment associated with an SoC range. 
\begin{figure}[t]
    \centering
    \subfloat[Marginal SoC opportunity value and bids.]{
    \includegraphics[trim = 5mm 0mm 0mm 0mm, clip, width = .9\columnwidth]{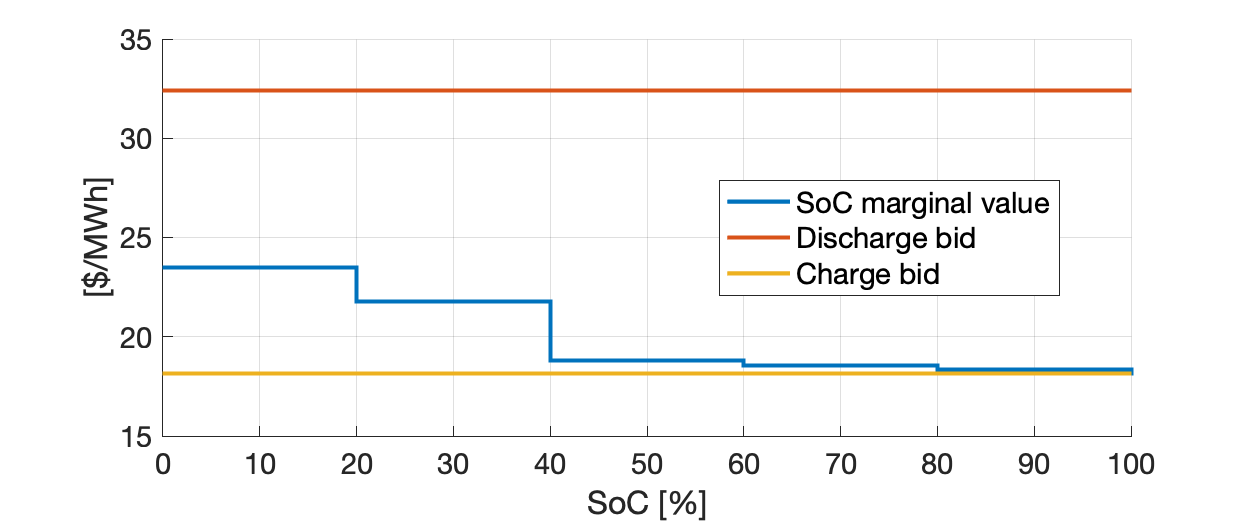}
    \label{fig:soc_bid1}
    }
    \\
    \subfloat[SoC opportunity value.]{
    \includegraphics[trim = 5mm 0mm 0mm 0mm, clip, width = .9\columnwidth]{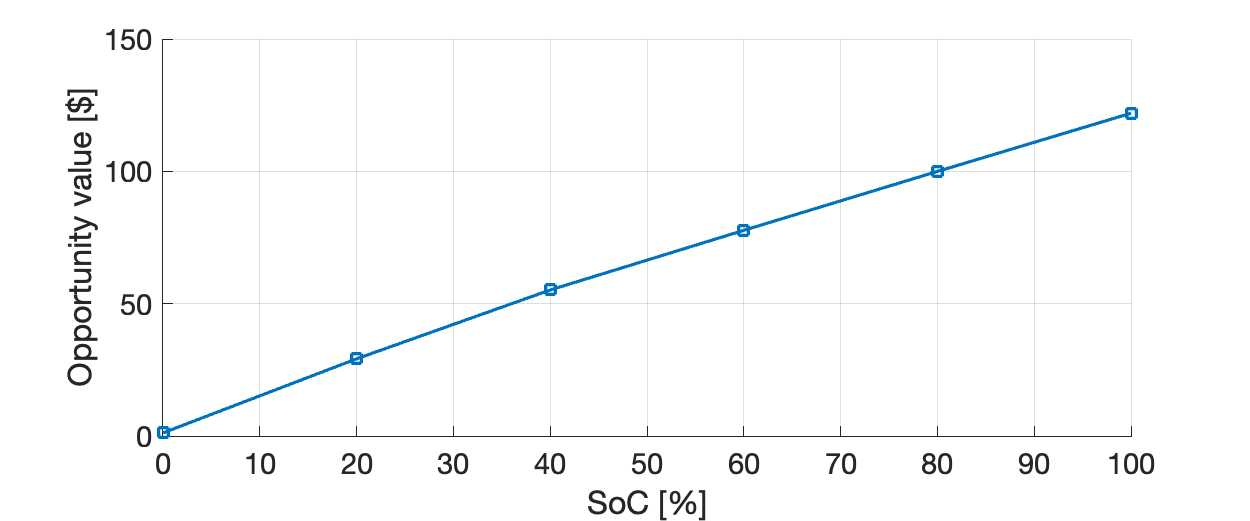}
    \label{fig:soc_bid2}
    }
    \caption{A 5-segment power and SoC bid example of a 1~MW/6~MWh storage. The discharge cost is \$10/MWh and the round-trip efficiency is 81\%.}
    \label{fig:soc_bid}
\end{figure}

In the SoC bid model, a storage participant  $m$ submits bids  $c\up{e}_{m,j,t}$ with $j\in J$ segments representing the marginal opportunity value associated with each SoC range $[E_{m,j-1}, E_{m,j}]$. $E_{m,0}$ represents zero SoC (or the minimum allowed SoC) and $E_{m,J}$ represents full SoC (or the maximum allowed SoC). In addition, storage participants must submit a discharge cost $c_m$ which covers the cost of degradation and other O\&M cost, discharge efficiency $\eta$. The economic dispatch model with storage SoC bids is formulated as
\begin{subequations}\label{ed2}
\begin{align}
    &\min_{p_{m,t}, b_{m,t}, g_{k,t}, e_{m,t}} \sum_{k} C_k(g_{k,t})  + \sum_{m} c_m p_{m,t} - v_{m,t} \label{ed2:obj}\\
    &\text{subjects to }\nonumber\\
    &0 \leq p_{m,t} \leq P\up{p}_m(e_{m,t-1}) \label{ed2:c1}\\
    &0 \leq b_{m,t} \leq P\up{b}_m(e_{m,t-1}) \label{ed2:c2}\\
    & e_{m,t} = e_{m,t-1} - p_{m,t}/\eta_m + b_{m,t}\eta_m \label{ed2:c3}\\
        & E_{m,0} \leq e_{m,t} \leq E_{m,J} \label{ed2:c5}\\
    & v_{m,t} \leq c\up{e}_{m,j,t}(e_t-E_j) + \sum_{\tau=1}^{j-1} c\up{e}_{m,\tau,t}(E_{m,\tau}-E_{m,\tau-1})  \label{ed2:c4}\\
    & \text{other generator, network, nodal balance constraints}\nonumber
\end{align}
\end{subequations}
in which the cost of storage dispatch is represented using the term $c_m p_{m,t} - v_{m,t}$ for each storage $m$. $c_m$ is the storage discharge cost, $v_{m,t}$ models the change in storage opportunity value, it has a minus sign in a minimization optimization representing the system operator trade-off between reducing current system operating cost and increasing the opportunity value of energy stored in the storage. \eqref{ed2:c1} and \eqref{ed2:c2} model the charge and discharge power rating, same as in the power bids model. \eqref{ed2:c3} models the change in the storage state-of-charge during the dispatch step, again, note that since this is a single-period dispatch problem, 
$e_{m,t-1}$ is an input parameter to the economic dispatch representing the monitored storage SoC at the beginning of the dispatch period. \eqref{ed2:c5} enforces the upper and lower SoC limits for storage $m$. \eqref{ed2:c4} models the piece-wise linear storage opportunity value function, the first term on the right-hand side models the starting opportunity value of the storage at the beginning of segment $j$, while the second term models the increasing slope of the opportunity value between the segment. Fig.~\ref{fig:soc_bid} shows an example of a 5-segment SoC opportunity bid curve, note that the marginal opportunity value decreases monotonically with the SoC, hence the opportunity value function (the integral of the marginal opportunity value) as shown in Fig.~\ref{fig:soc_bid2} is a concave function.

SoC bids $c\up{e}_{m,j,t}$ can be directly derived using conditional expectation of the marginal value function $q_t(e)$ based on the SoC segment ranges
\begin{align}
    c\up{e}_{m,j,t} = \frac{1}{E_{m,j}-E_{m,j-1}}\int_{E_{m,j-1}}^{E_{m,j}}q_t(e) de
\end{align}
note that we can also generalize the power bid model as an SoC bid model with only one SoC segment, in which case $c\up{e}_{m,j,t}$ becomes the average of the marginal opportunity value. 

\section{Case Study}
We test combinations of different storage valuation methods and bidding models with historical data from New York Independent System Operator (NYISO). We assume energy storage are price takers, and use historical price data to simulate storage dispatch.

\subsection{Arbitrage-based Economic Dispatch Simulation}
We adopt price-taker assumptions for storage participating in real-time dispatch. We assume storage participants do not seek to exercise market power and have negligible influences on market clearing prices. Under price-taker assumptions, dispatch results from profit-maximizing price arbitrage are equivalent to cost-minimizing economic dispatch results, and profit collected by the storage also reflects their utilization in reducing total system operating costs~\cite{castillo2013profit}. Therefore, we will use price arbitrage to simulate how would storage perform in historical system operations, and use the arbitrage profit as an indicator of storage utilization in the power system. The higher the profit, the better utilization of the storage.

We define historical real-time price data as $\pi_t$, different from the price $\lambda_t$ used in storage opportunity valuation as in \eqref{eq1}. Here,  $\lambda_t$ represents the prediction of $\pi_t$, hence we use price predictions to perform storage valuation, and use the actual price to simulate storage arbitrage. Under the price-taker assumptions, storage dispatch from the economic dispatch model with power bids \eqref{ed1} is equivalent to results from the following single-period arbitrage model:
\begin{align}\label{arb1}
    \max_{p_{m,t}, b_{m,t}} \pi_t(p_{m,t}-b_{m,t}) - c\up{p}_{m,t} p_{m,t} + c\up{b}_{m,t} b_{m,t}
\end{align}
subject to constraints \eqref{ed1:c1} and \eqref{ed1:c2}.

While the equivalent arbitrage model for the economic dispatch model with SoC bids \eqref{ed2} is 
\begin{align}\label{arb2}
    \max_{p_{m,t}, b_{m,t}} \pi_t(p_{m,t}-b_{m,t}) - c\up{p} p_{m,t} + v_{m,t}
\end{align}
subject to constraints \eqref{ed2:c1}--\eqref{ed2:c4}.

Note that we only need to consider one storage in the arbitrage optimization because of the price-taker assumption, that storage will not affect market clearing prices, so different storage will not affect each other. We perform arbitrage optimizations in sequential order and update the storage SoC after each time period.

\subsection{Price Data}

\begin{figure}[t]%
	\centering
	\subfloat[New York State Price Zones.]{
		\includegraphics[trim = 0mm 0mm 0mm 0mm, clip, width = .85\columnwidth]{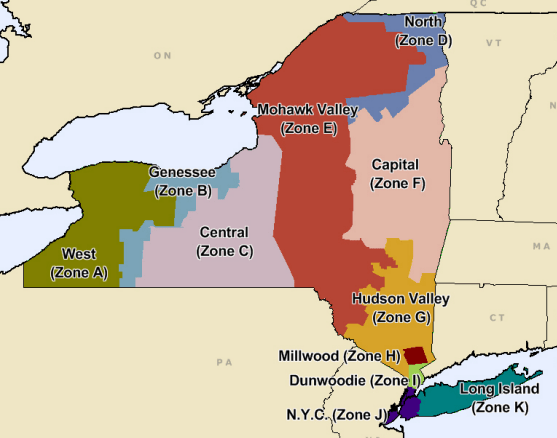}
		\label{fig:nyc1}%
	}
	\\
	\subfloat[30 day moving average price.]{
		\includegraphics[trim = 5mm 0mm 10mm 0mm, clip, width = .95\columnwidth]{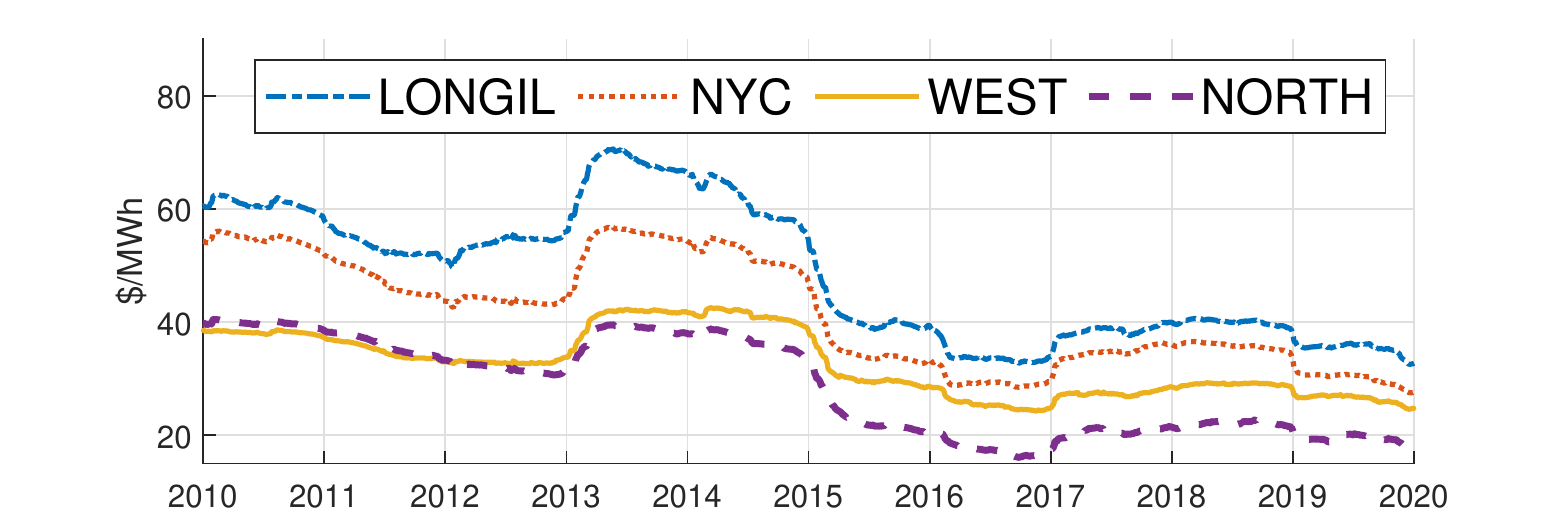}
		\label{fig:nyc2}%
	}
	\\
	\subfloat[30 day moving average daily price deviations.]{
		\includegraphics[trim = 5mm 0mm 10mm 0mm, clip, width = .95\columnwidth]{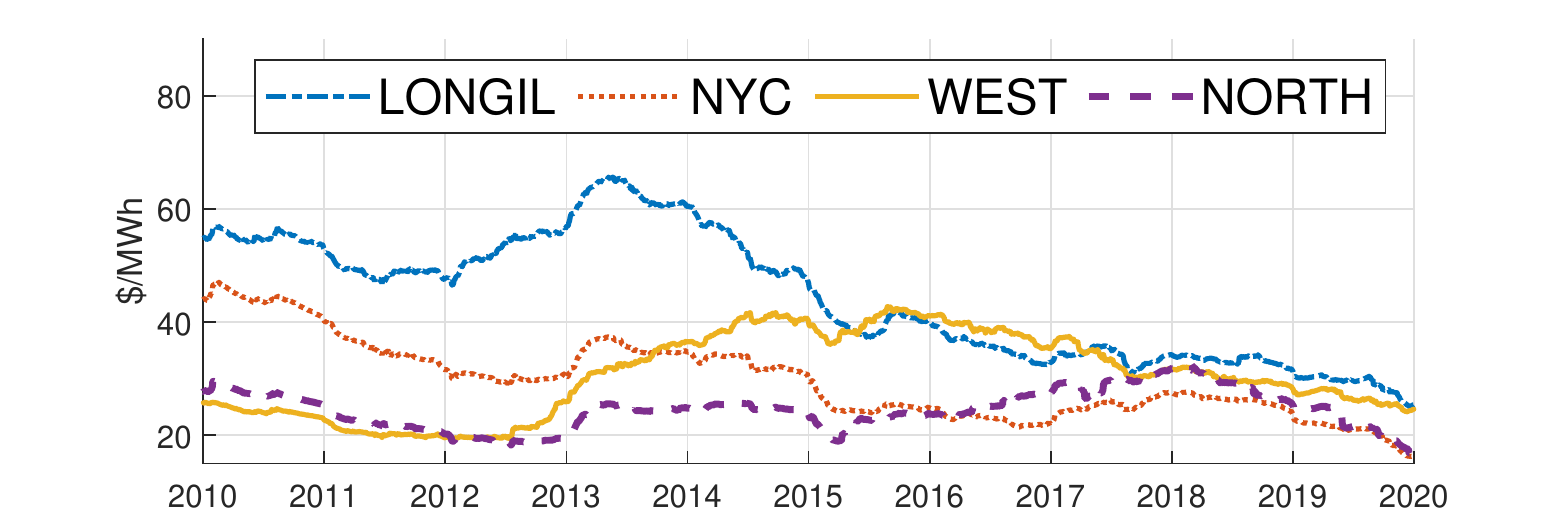}
		\label{fig:nyc3}%
	}
  \caption{Location of price zones and the historical price data of the selected four zones.}
    \label{fig:nyc}
\end{figure}

We use historical price data in 2019 from the WEST, NORTH, NYC, and LONGIL zones in NYISO. Fig.~\ref{fig:nyc} shows the geographical location of these four zones and overviews of historical price mean and standard deviations. These four zones have distinct generation resource mix in 2019: WEST was almost 100\% hydro, NORTH was about 10\% wind and 90\% hydro, generations in NYC and LONGIL were dominated by natural gas, while LONGIL has a higher ratio of oil generations and also suffers severe congestion~\cite{nyiso}. As a result, LONGIL has the highest average price and price deviations, WEST and NORTH also have high price deviations likely due to high renewable share, but the average price is low. NYC has the second high average price, but the price deviation is the lowest. 

\begin{figure*}[t]%
	\centering
	\subfloat[Day-ahead forecast bids, 1 hour storage]{
		\includegraphics[trim = 5mm 0mm 10mm 0mm, clip, width = .95\columnwidth]{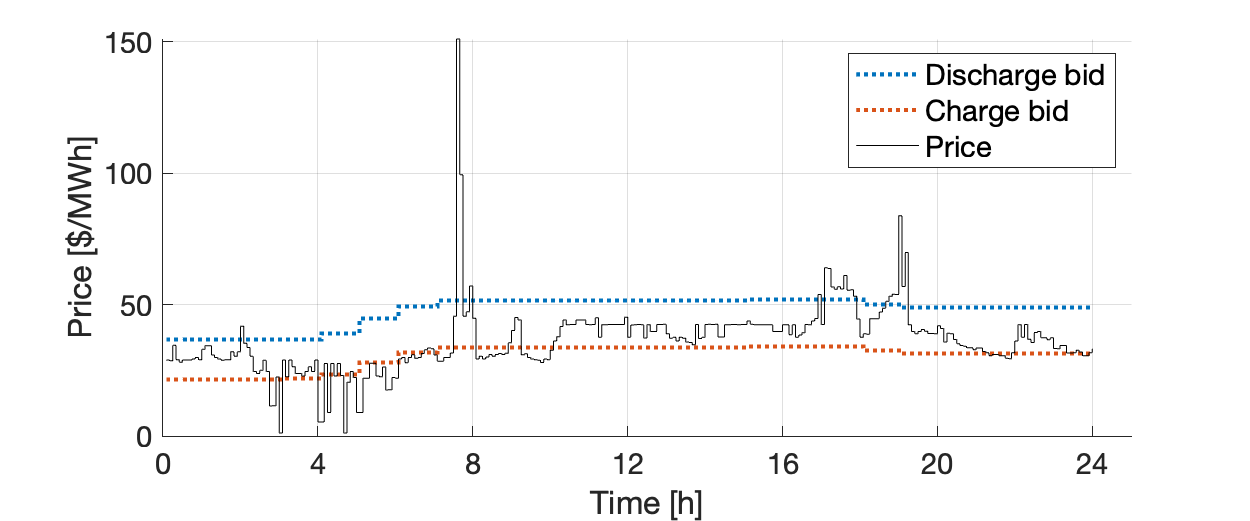}
		\label{fig:bid1}%
	}
	\subfloat[Perfect forecast bids, 1 hour storage]{
		\includegraphics[trim = 5mm 0mm 10mm 0mm, clip, width = .95\columnwidth]{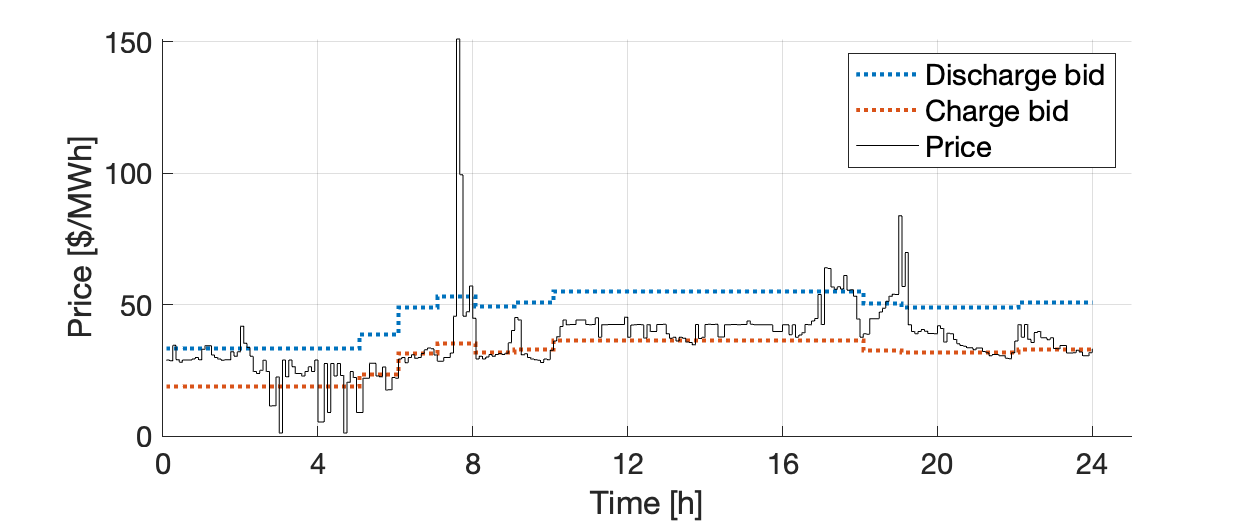}
		\label{fig:bid2}%
	}
	\\
	\subfloat[Day-ahead forecast bids, 6 hour storage]{
		\includegraphics[trim = 5mm 0mm 10mm 0mm, clip, width = .95\columnwidth]{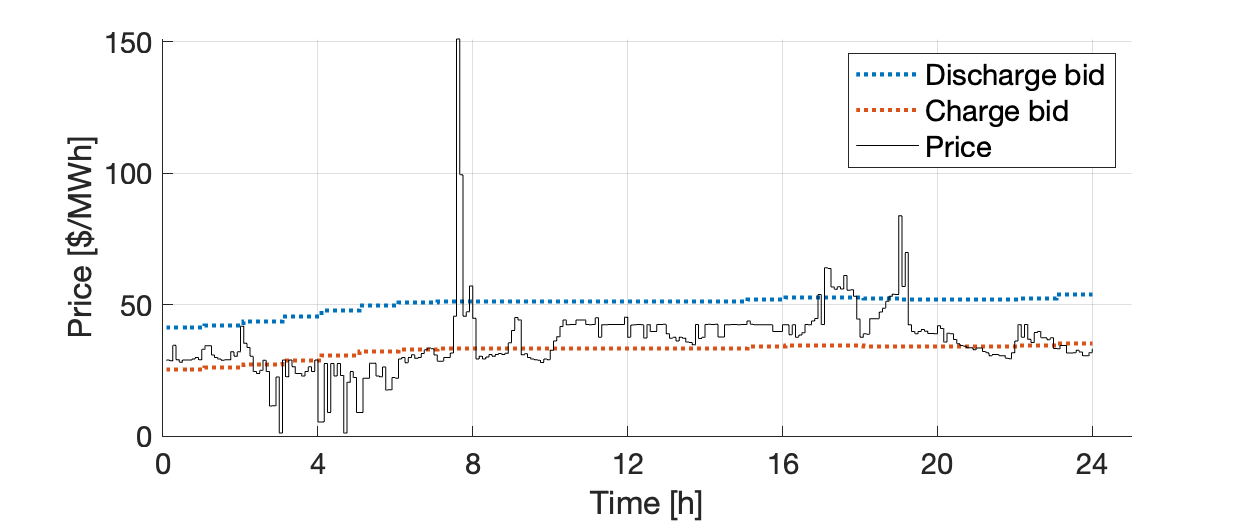}
		\label{fig:bid3}%
	}
	\subfloat[Perfect forecast bids, 6 hour storage]{
		\includegraphics[trim = 5mm 0mm 10mm 0mm, clip, width = .95\columnwidth]{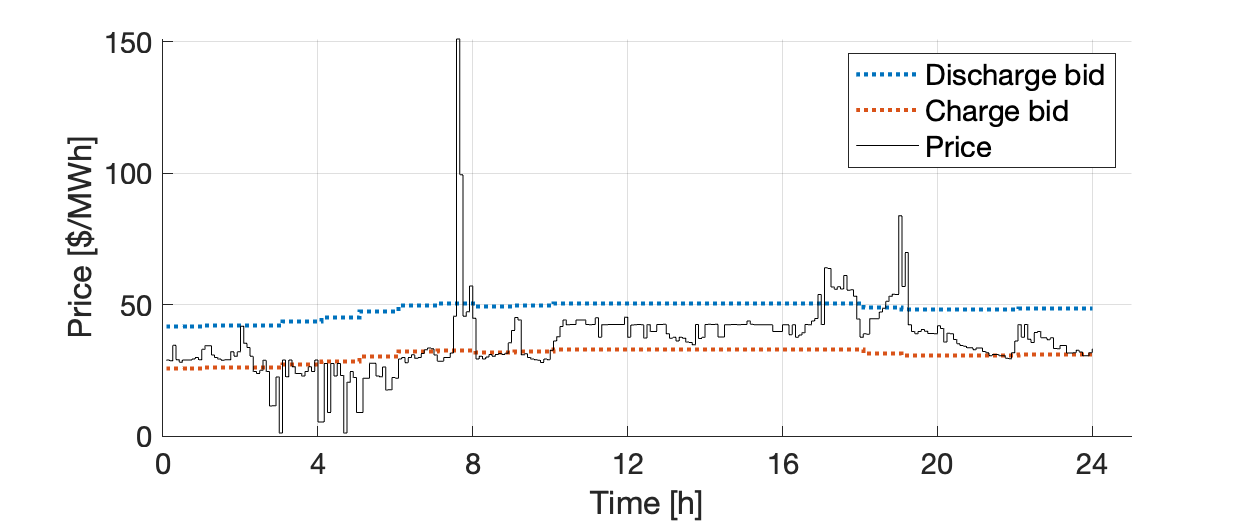}
		\label{fig:bid4}%
	}
	\\
	\subfloat[Day-ahead forecast bids, 24 hour storage]{
		\includegraphics[trim = 5mm 0mm 10mm 0mm, clip, width = .95\columnwidth]{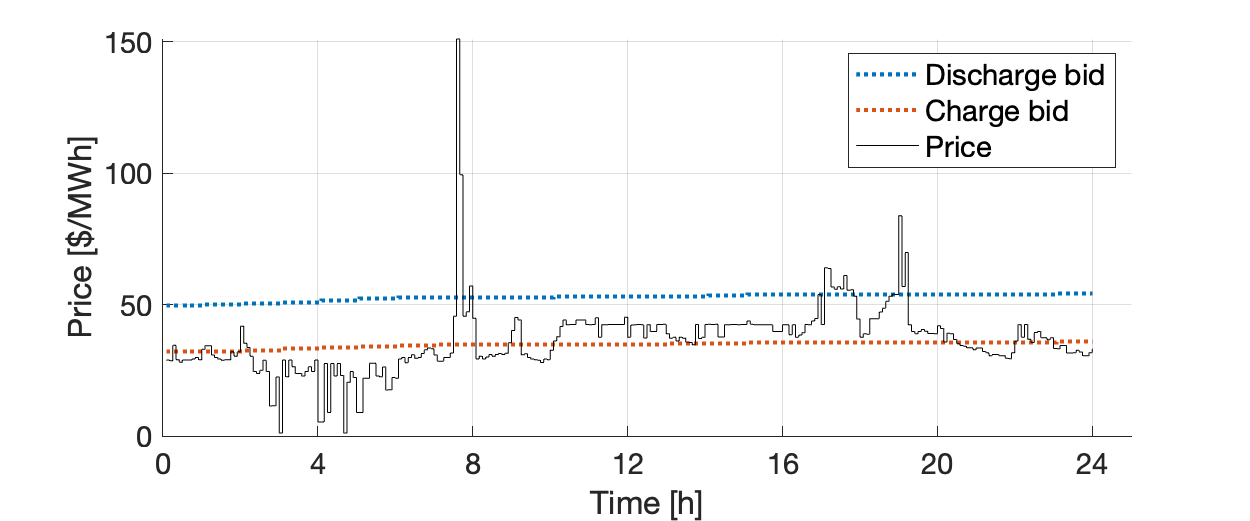}
		\label{fig:bid5}%
	}
	\subfloat[Perfect forecast bids, 24 hour storage]{
		\includegraphics[trim = 5mm 0mm 10mm 0mm, clip, width = .95\columnwidth]{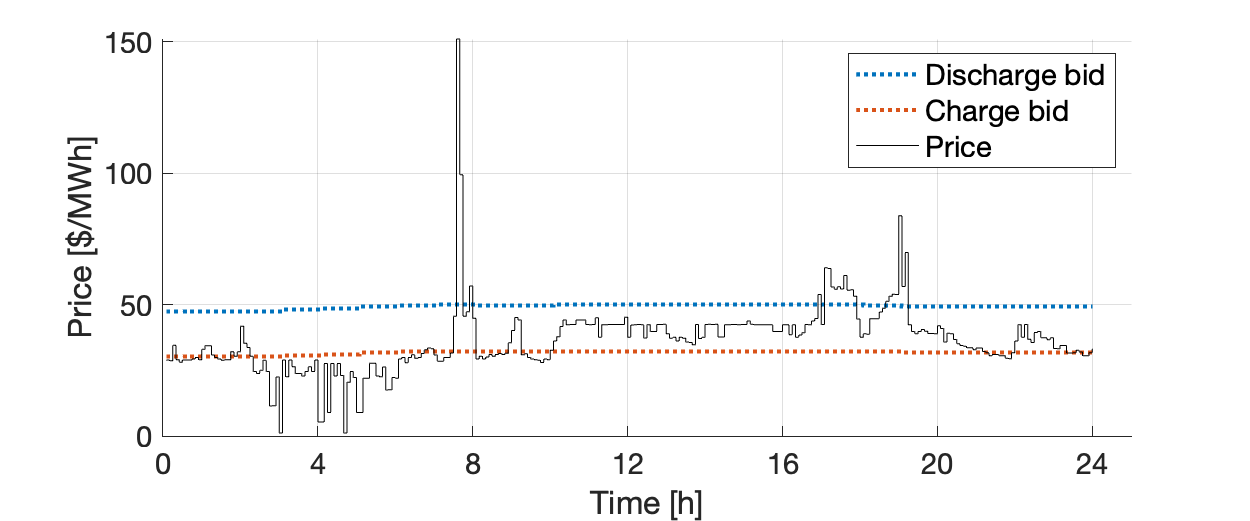}
		\label{fig:bid6}%
	}
  \caption{Comparison of storage charge and discharge bids from different price predictions and storage durations.}
    \label{fig:bid}
\end{figure*}

\begin{figure*}[t]%
	\centering
	\subfloat[LONGIL discharge bids]{
		\includegraphics[trim = 0mm 0mm 5mm 10mm, clip, width = .48\columnwidth]{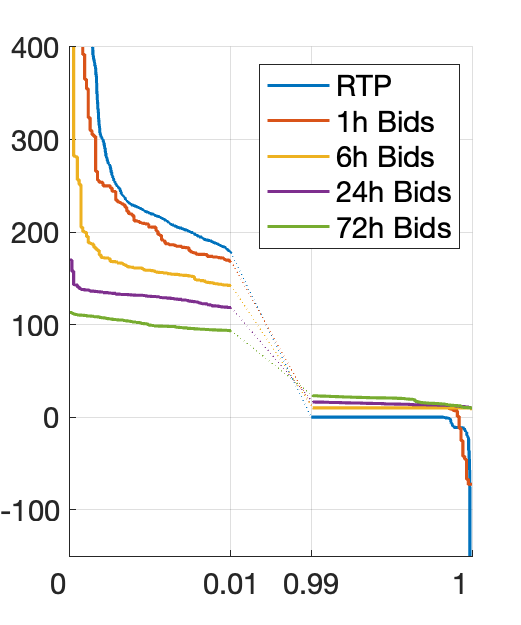}
		\label{fig:dur1}%
	}
	\subfloat[NORTH discharge bids]{
		\includegraphics[trim = 0mm 0mm 5mm 10mm, clip, width = .48\columnwidth]{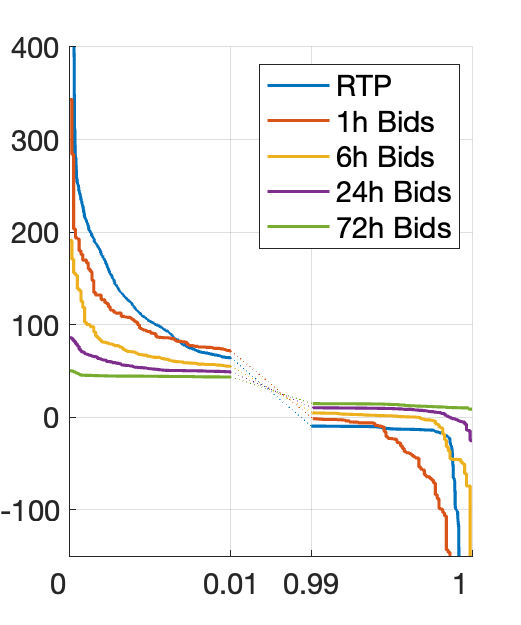}
		\label{fig:dur2}%
	}
	\subfloat[NYC discharge bids]{
		\includegraphics[trim = 0mm 0mm 5mm 10mm, clip, width = .48\columnwidth]{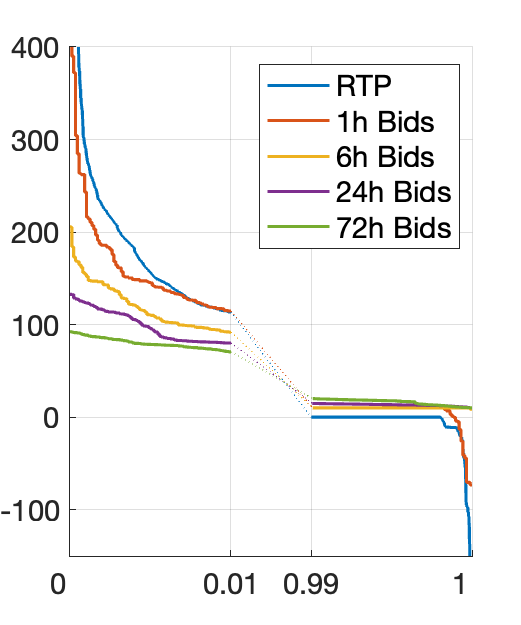}
		\label{fig:dur3}%
	}
	\subfloat[WEST discharge bids]{
		\includegraphics[trim = 0mm 0mm 5mm 10mm, clip, width = .48\columnwidth]{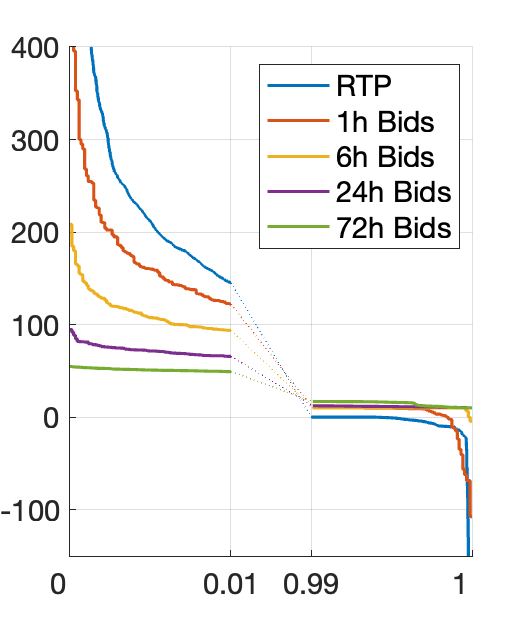}
		\label{fig:dur4}%
	}
	\\
	\subfloat[LONGIL charge bids]{
		\includegraphics[trim = 0mm 0mm 5mm 10mm, clip, width = .48\columnwidth]{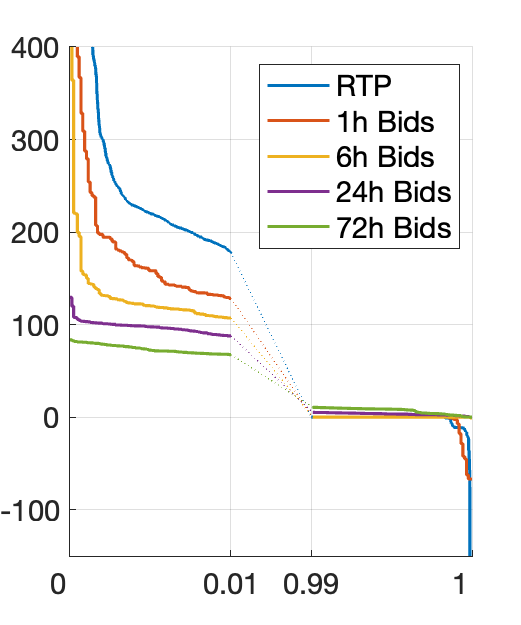}
		\label{fig:dur5}%
	}
	\subfloat[NORTH charge bids]{
		\includegraphics[trim = 0mm 0mm 5mm 10mm, clip, width = .48\columnwidth]{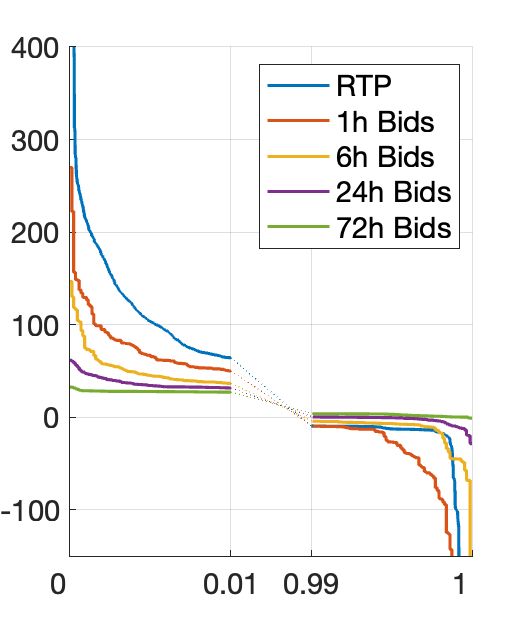}
		\label{fig:dur6}%
	}
	\subfloat[NYC charge bids]{
		\includegraphics[trim = 0mm 0mm 5mm 10mm, clip, width = .48\columnwidth]{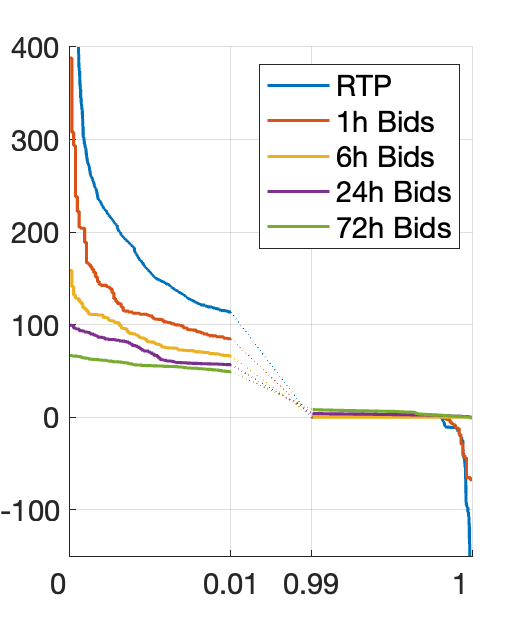}
		\label{fig:dur7}%
	}
	\subfloat[WEST charge bids]{
		\includegraphics[trim = 0mm 0mm 5mm 10mm, clip, width = .48\columnwidth]{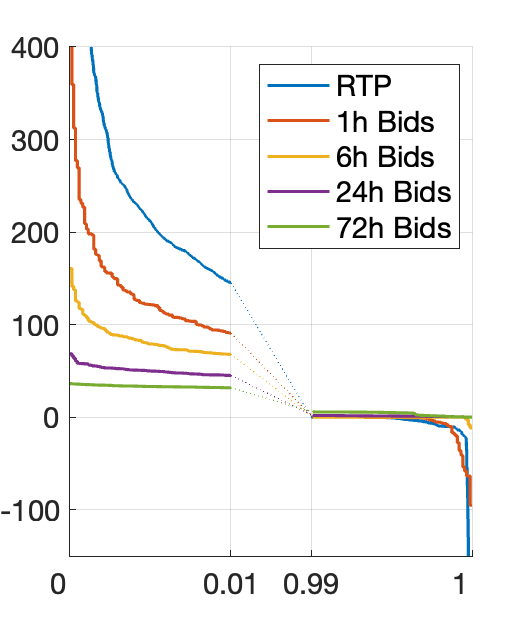}
		\label{fig:dur8}%
	}
	\\
	\subfloat[LONGIL average bids]{
		\includegraphics[trim = 0mm 0mm 5mm 10mm, clip, width = .48\columnwidth]{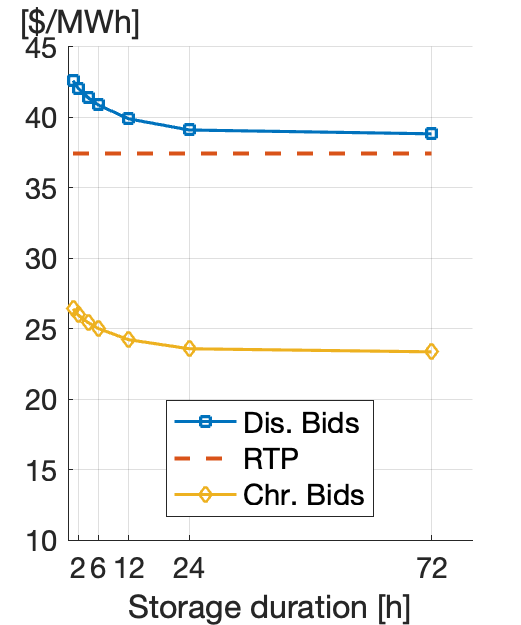}
		\label{fig:dur9}%
	}
	\subfloat[NORTH average bids]{
		\includegraphics[trim = 0mm 0mm 5mm 10mm, clip, width = .48\columnwidth]{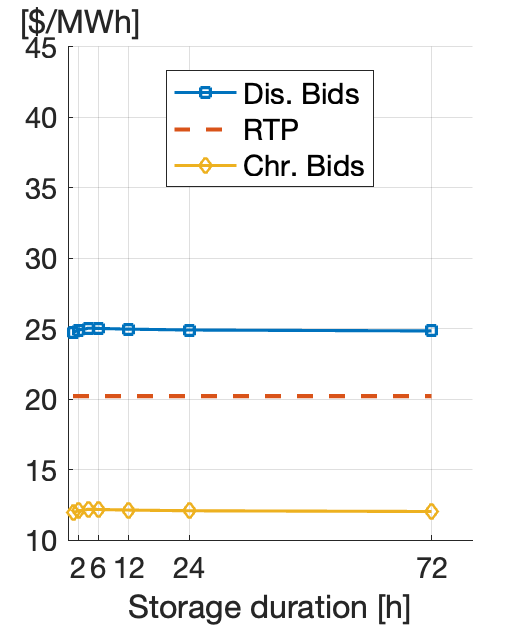}
		\label{fig:dur10}%
	}
	\subfloat[NYC average bids]{
		\includegraphics[trim = 0mm 0mm 5mm 10mm, clip, width = .48\columnwidth]{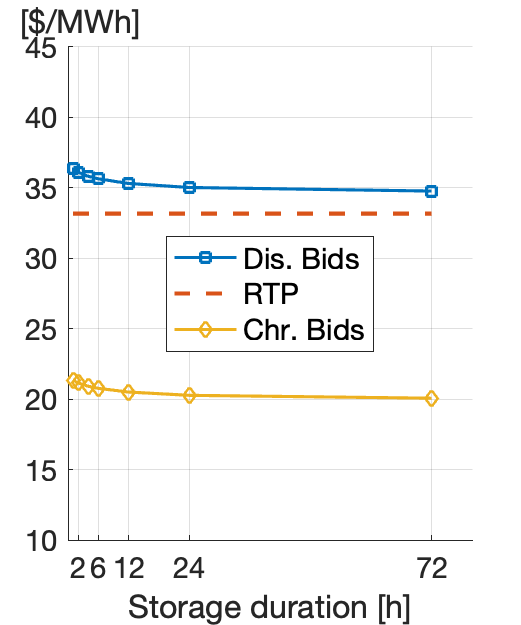}
		\label{fig:dur11}%
	}
	\subfloat[WEST average bids]{
		\includegraphics[trim = 0mm 0mm 5mm 10mm, clip, width = .48\columnwidth]{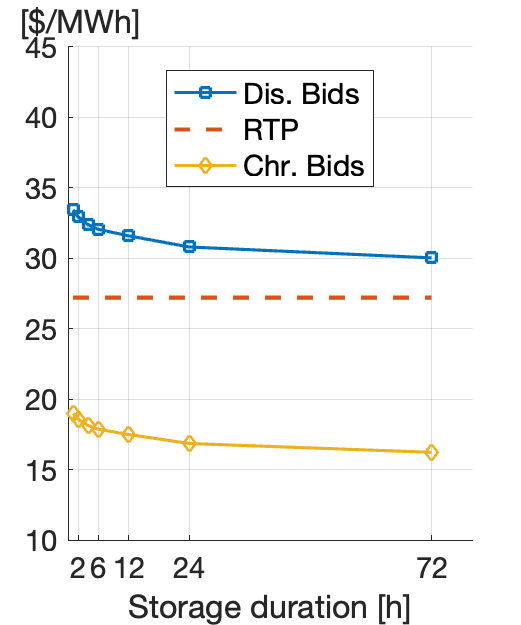}
		\label{fig:dur12}%
	}
  \caption{Price duration curves, average of real-time price (RTP) and storage discharge (Dis.) and charge (Chr.) bids from different locations and storage duration. Bids are generated using perfect real-time price predictions. The duration curves have been modified to highlight the top and last 1\% quantile of the distribution. Units in all y axes are \$/MWh.}
    \label{fig:dur}
\end{figure*}

\subsection{Dispatch Test Setting}
We design the following test cases to compare how utilization of storage varies with different market designs and valuation methods (price forecasts):
\begin{enumerate}
    \item \textbf{DA-PB-DF}: The storage uses day-ahead price (DF) to design power bids (PB) and participate in day-ahead (DA) markets. We use day-ahead price in \eqref{eq1} for valuation and use day-ahead price in the power bid arbitrage model \eqref{arb1}.
    \item \textbf{DA-SB-DF}: The storage uses day-ahead price (DF) to design SoC bids (SB) and participate in day-ahead (DA) markets. We use day-ahead price in \eqref{eq1} for valuation and use day-ahead price in the power bid arbitrage model \eqref{arb2}. 
    \item \textbf{RT-PB-DF}: The storage uses day-ahead price as real-time price predictions (DF) to design power bids (PB) and participate in real-time (RT) markets. We use day-ahead price in \eqref{eq1} for valuation and use real-time price in the power bid arbitrage model \eqref{arb1}.
    \item \textbf{RT-SB-DF}: The storage uses day-ahead price as real-time price predictions (DF) to design SoC bids (SB) and participate in real-time (RT) markets. We use day-ahead price in \eqref{eq1} for valuation and use real-time price in the SoC bid arbitrage model \eqref{arb2}.
    \item \textbf{RT-PB-PF}: The storage uses perfect real-time price forecasts (PF) to design power bids (PB) and participate in real-time (RT) markets. We use real-time price in \eqref{eq1} for valuation and use real-time price in the power bid arbitrage model \eqref{arb1}.
    \item \textbf{RT-SB-PF}: The storage uses perfect real-time price forecasts (PF) to design SoC bids (SB) and participate in real-time (RT) markets. We use real-time price in \eqref{eq1} for valuation and use real-time price in the SoC bid arbitrage model \eqref{arb2}. Note that when using the same price data to perform valuation and to execute arbitrage using SoC bids,  it is equivalent to solving a multi-period arbitrage optimization using dynamic programming, hence the result is the same as directly solving a multi-period arbitrage with perfect price data. Hence, arbitrage revenue from this case represents the highest utilization possible for storage.
\end{enumerate}
Day-ahead prices are of hourly resolution and real-time prices are of 5-minute resolution. All bids are of hourly resolution.
For SoC bids, we use 20 segments per hour of duration, corresponding to existing generator bidding rules that each hourly bid can contain 20 segments. 
We use day-ahead price as predictions of real-time prices. Although this is a naive approach for real-time price predictions, we use this as a lower bound for investigating storage participation models based on price predictions, as the day-ahead price is available to the public.  

\begin{figure*}[t]%
	\centering
	\subfloat[NYC]{
		\includegraphics[trim = 5mm 0mm 10mm 0mm, clip, width = .95\columnwidth]{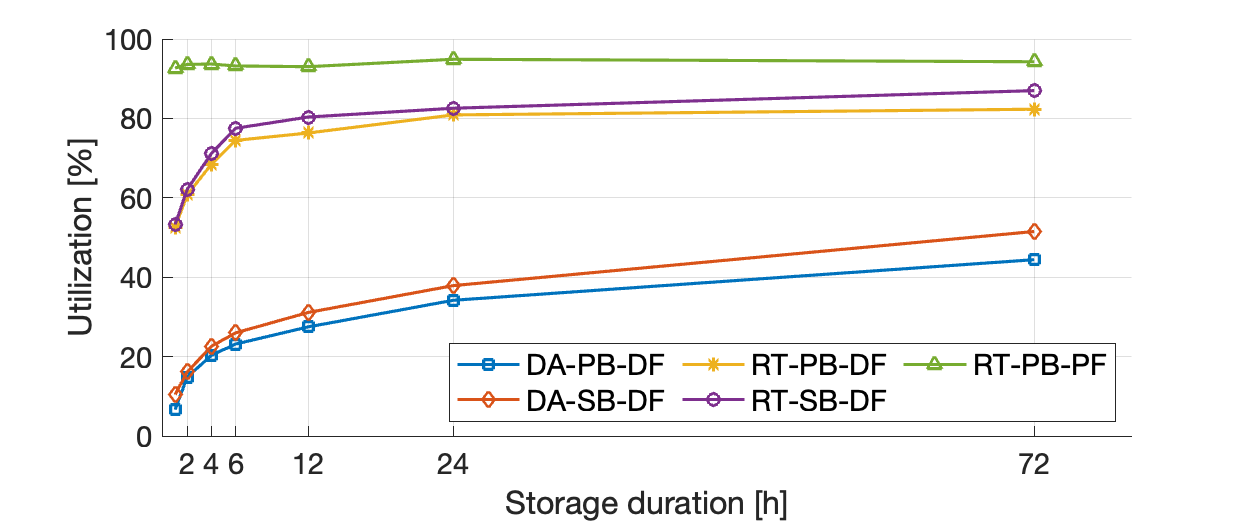}
		\label{fig:ut1}%
	}
	\subfloat[LONGIL]{
		\includegraphics[trim = 5mm 0mm 10mm 0mm, clip, width = .95\columnwidth]{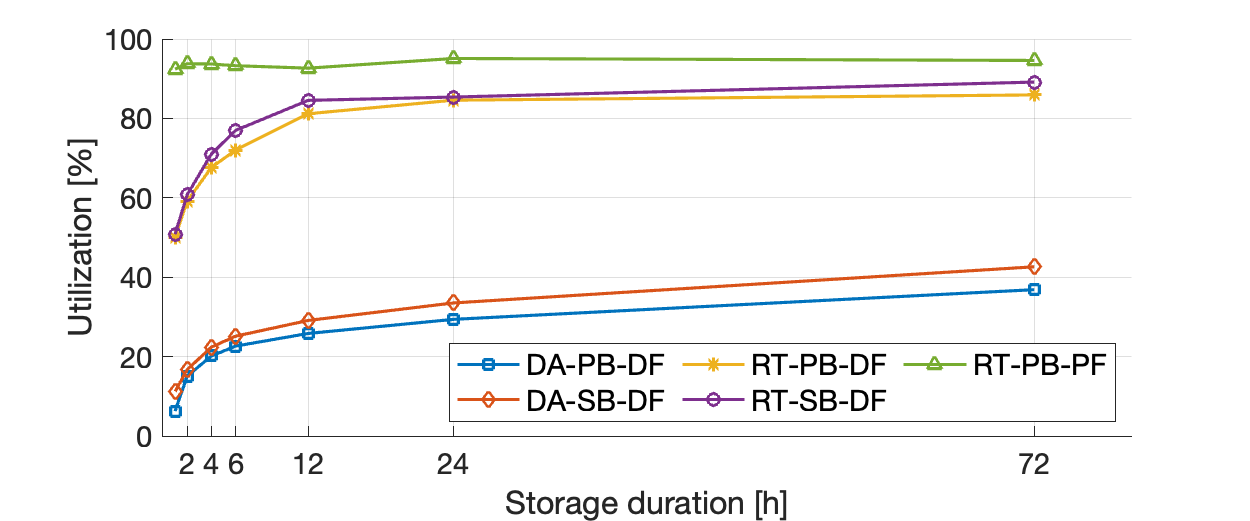}
		\label{fig:ut2}%
	}
	\\
	\subfloat[WEST]{
		\includegraphics[trim = 5mm 0mm 10mm 0mm, clip, width = .95\columnwidth]{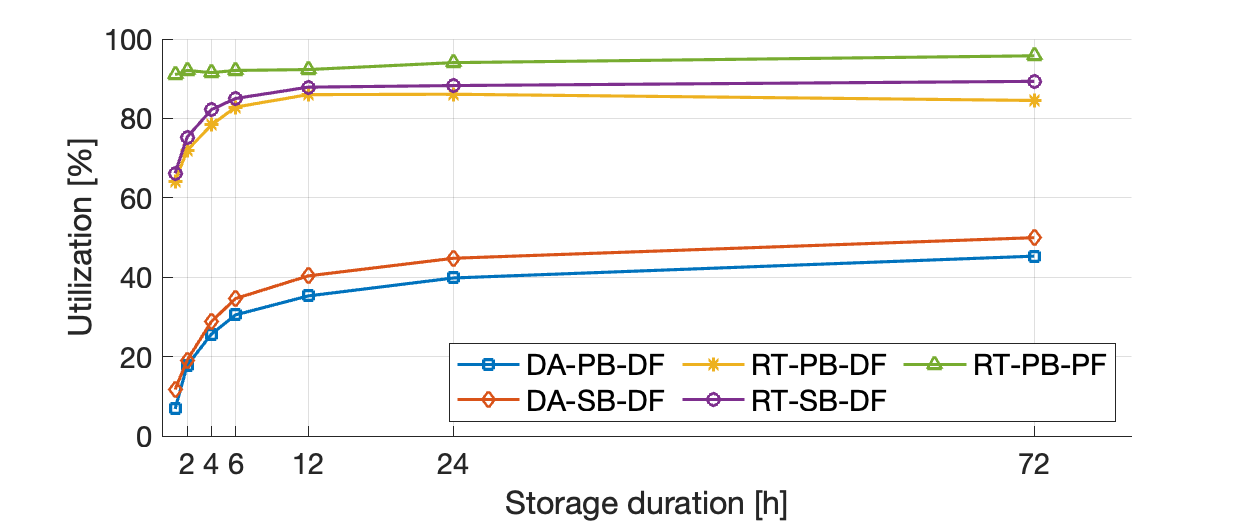}
		\label{fig:ut3}%
	}
	\subfloat[NORTH]{
		\includegraphics[trim = 5mm 0mm 10mm 0mm, clip, width = .95\columnwidth]{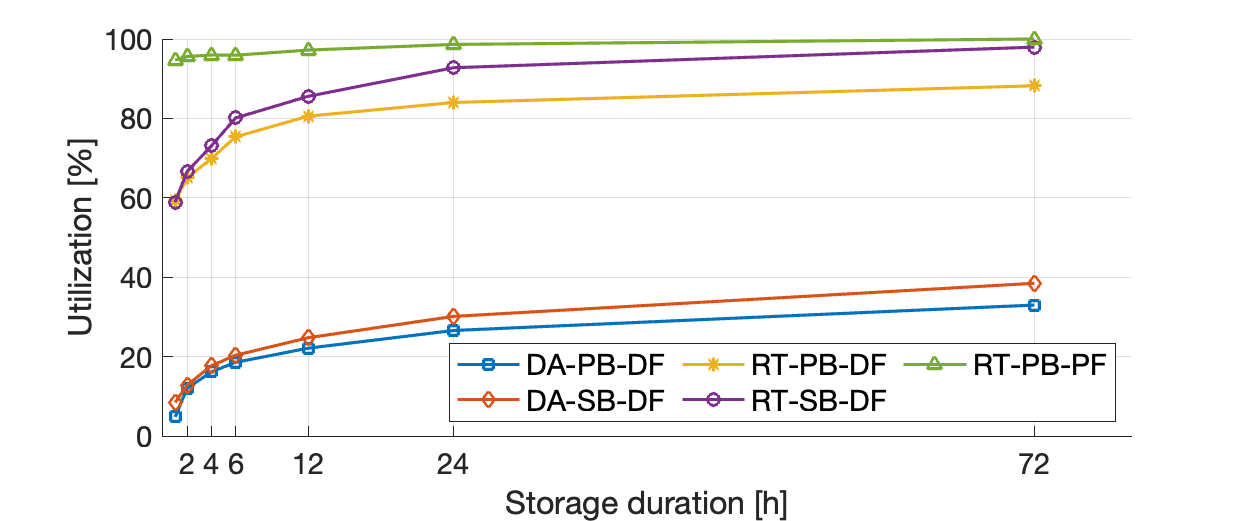}
		\label{fig:ut4}%
	}
  \caption{Comparison of storage utilization in four NYISO price zones in single-period economic dispatches, the utilization rate is with respect to multi-period dispatch with perfect forecasts.}
    \label{fig:ut}
\end{figure*}

\subsection{Storage Settings}
In the case studies, we consider general storage models with durations of one, two, four, six, twelve, twenty-four, and seventy-two hours of duration. We assume the storage has a round-trip efficiency of 81\% and has a fixed discharge cost of \$10/MWh, covering the cost of degradation and other O\&M costs. 

\section{Results}

\subsection{Impact of Forecast and Duration over Storage Bids}

We compare storage charge and discharge power bids in real-time market using day-ahead price forecast or perfect price forecast over different storage durations. Fig.~\ref{fig:bid} shows a comparison of bids from  1-hour, 6-hour, and 24-hour duration storage. Bids from the 1-hour storage are more sensitive to price volatilities, especially in the perfect price forecast case,  storage bids follow the general trend of the real-time price. Bids designed using day-ahead prices are less volatile because day-ahead prices have fewer fluctuations.

As storage capacity increases, variations in bids reduce. In the 6-hour duration storage case, daily bid variation - the difference between the daily maximum bid and the minimum bid values - is only about half compared to the 1-hour storage case. In the 24-hour storage duration case, bids become almost the same throughout the day. In all cases, the gap between the bids is stable as it covers storage discharge cost and efficiency losses.

Figure~\ref{fig:dur} shows the price duration curves of real-time prices and the power bids from different storage durations generated using the real-time price data. The results echo with Fig.~\ref{fig:bid} that bids become less volatile as the storage duration increases. The real-time price data from all four zones shows clear occurrences of price spikes and negative prices, with LONGIL having the highest spikes and NORTH having the most negative prices. Discharge bids from one-hour storage still follow price spike trends, especially in LONGIL that storage could still bid higher than \$400/MWh. However, after increasing storage duration to more than six hours, storage discharge bids become smoother, all bids are less than \$200/MWh except in LONGIL. For the 72-hour storage, discharge bids look almost like a linear at the 1\% and 99\% quantile, indicating bids from storage are very smooth.

Finally, we compare the annual average of discharge and charge bids and the real-time price, shown in Fig.~\ref{fig:dur9}--\ref{fig:dur12}. The gap between discharge and charge bids includes the \$10/MWh discharge cost and the 81\% round trip efficiency. Both charge and discharge bids values reduce as the storage capacity increase,  discharge bids are slightly higher than average real-time prices except in NORTH zone, which has a much lower average price and more occurrence of negative prices due to the higher share of wind generations. Notably, the higher the average real-time price, the closer the discharge bids to the real-time price. These results, together with the price duration curve results, show that the power bids, designed using a price-taker arbitrage profit-maximization model, will help to mitigate extreme prices but likely keep the average price similar to without storage.


\subsection{Impact of Forecast and Duration over Storage Utilization}

Fig.~\ref{fig:ut} shows the comparison of storage utilization rate in the four NYISO price zones with different market design and participation settings. We normalized the profit from other cases with the profit from the perfect price forecast and dispatch (RT-SB-PF) case, which represents the highest amount of profit possible to obtain. Hence, the percentage value shown in these figures represents how much utilization can be achieved compared to the maximum possible. 

The first observation is that utilization of storage in day-ahead markets is low over all price zones, less than 40\% in most cases. This is not surprising as day-ahead prices usually have much fewer volatilities compared to real-time prices. The value of the storage is not on generating electricity, but on shifting the energy. Hence, dispatch storage in day-ahead markets will achieve much less utilization of the storage resource compared to real-time markets. However, we must note that this conclusion is drawn based on the price-taker assumptions that storage has a negligible impact on market clearing results. Because day-ahead market clearing includes unit commitments, dispatch storage in the day-ahead market may change the commitment decision over some thermal generations, letting them turn off early or not be turned on at all. This would provide higher cost savings to the system operation but may not reflect in the price, as generator start-up costs and no-load costs are usually not included in market clearing prices~\cite{kirschen2018fundamentals}. 

For the two bidding cases using day-ahead prices as real-time price predictions (RT-PB-DF and RT-SB-DF), the utilization rate is low with short storage durations - between 50\% to 60\% for 1-hour storage, but increases quickly as the storage duration increases. In all price zones, the utilization rate reached above 80\% when the storage duration exceeds 12-hour, and settled around 90\% in the 72-hour duration case. Recall that we are using day-ahead price as real-time price predictions which is a naive prediction approach, while the utilization rate could be further improved with more advanced price predictions or valuation approaches~\cite{zheng2022arbitraing}. Still, this naive approach achieves good performance with long-duration storage, in particular, in the NORTH zone, the 72-hour storage reached almost 100\% utilization rate with the SoC bid dispatch model. 

\begin{figure}[t]
    \centering
    \includegraphics[trim = 5mm 0mm 10mm 0mm, clip, width = .95\columnwidth]{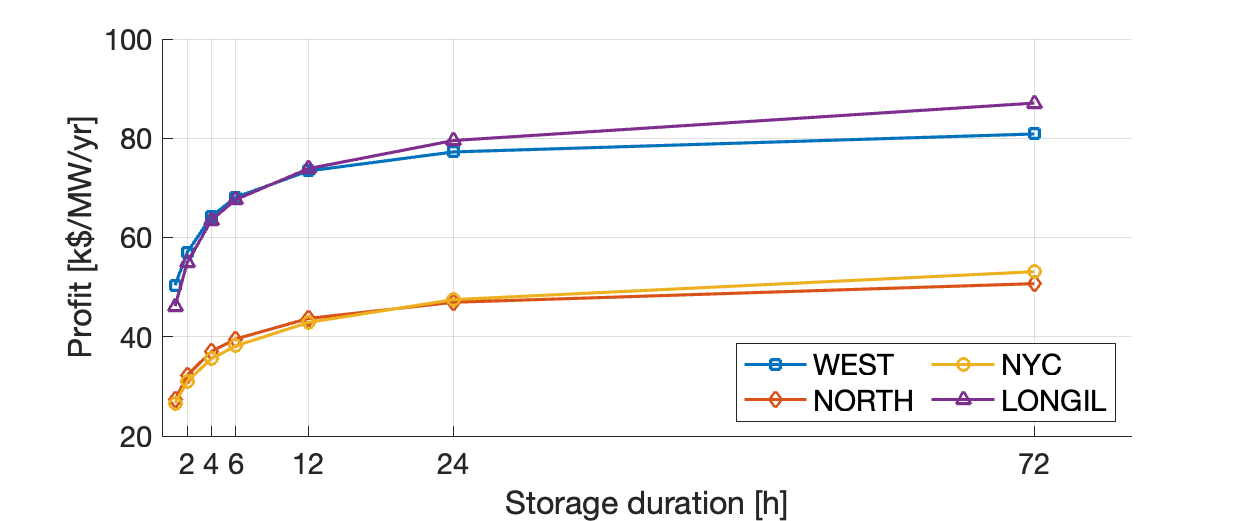}
    \caption{Comparison of perfect price prediction arbitrage profit across four NYISO price zones.}
    \label{fig:pro}
\end{figure}

We compare the difference in storage utilization between the power bidding model and the SoC bidding model. Recall that in the power bidding model, the storage submits a single bid price for discharging and a single bid price for charging, the system operator does not model storage SoC in the real-time dispatch. In the SoC bidding model, the storage bid is dependent on the storage SoC, and the system operator models the change in the SoC over the single dispatch time period. From the results, using storage bidding models provides around 5\% utilization improvement in all cases,  including day-ahead arbitrage using day-ahead forecast(DA-PB-DF), real-time arbitrage with day-ahead price forecasts (RT-PB-DF), real-time arbitrage using perfect forecasts (RT-PB-PF). Note that the utilization rate is with respect to the perfect real-time arbitrage results, if compare the power bids in day-ahead market (DA-PB-DF) with SoC bids in day-ahead markets (DA-SB-DF), the utilization rate is only around 80\%.

Finally, we compare utilization results across four different price zones. Storage utilization rates are the highest in NORTH zone compared to other zones, which is likely due to the more frequent occurrence of negative prices as the north zone has a higher wind generation ratio. Negative prices provide good opportunities as the storage can recharge while earning profit. Yet, comparing with the absolute profit from all four zones, shown in Fig.~\ref{fig:pro}, LONGIL and WEST provide nearly twice the profit compared to NYC and NORTH zones. Hence, storage will achieve better absolute utilization in these two zones. 

\section{Conclusion}

We used a price-taker arbitrage profit-maximization model to design energy storage bids into electricity markets, and investigated how different bidding and discharge models and storage durations impact the storage utilization rate. The duration curve and average value analysis over discharge and charge bids show that storage participants will help to mitigate extreme prices in electricity markets while the average price will likely remain similar to without storage, while storage with a longer duration will better contribute to flattening price fluctuations. 

We proposed an SoC bidding model in which the storage discharge and charge bids are dependent on SoC using a piece-wise linearization approach. The combination of perfect price prediction and the SoC bidding model achieves 100\% utilization of storage in single-period real-time dispatch with respect to dispatch storage using a deterministic multi-period optimization. Using SoC-independent charge and discharge bids will reduce the utilization rate by around 5\% regardless of the prediction accuracy, while using a naive day-ahead price prediction to design bids in real-time dispatch reduces the utilization rate to around 60\% for the one-hour storage duration, but utilization rate increases to near 90\% as the storage duration increases to 72 hours. Our result shows that price prediction accuracy is critical to integrating short-duration storage in real-time dispatch with less than four hours capacity. While long-duration storage with more than 12 hours of storage capacity can achieve a high utilization rate even with naive day-ahead price predictions. 


\bibliographystyle{IEEEtran}	
\bibliography{IEEEabrv,literature}		

\begin{thebibliography}{10}
\providecommand{\url}[1]{#1}
\csname url@samestyle\endcsname
\providecommand{\newblock}{\relax}
\providecommand{\bibinfo}[2]{#2}
\providecommand{\BIBentrySTDinterwordspacing}{\spaceskip=0pt\relax}
\providecommand{\BIBentryALTinterwordstretchfactor}{4}
\providecommand{\BIBentryALTinterwordspacing}{\spaceskip=\fontdimen2\font plus
\BIBentryALTinterwordstretchfactor\fontdimen3\font minus
  \fontdimen4\font\relax}
\providecommand{\BIBforeignlanguage}[2]{{%
\expandafter\ifx\csname l@#1\endcsname\relax
\typeout{** WARNING: IEEEtran.bst: No hyphenation pattern has been}%
\typeout{** loaded for the language `#1'. Using the pattern for}%
\typeout{** the default language instead.}%
\else
\language=\csname l@#1\endcsname
\fi
#2}}
\providecommand{\BIBdecl}{\relax}
\BIBdecl

\bibitem{shaner2018geophysical}
M.~R. Shaner, S.~J. Davis, N.~S. Lewis, and K.~Caldeira, ``Geophysical
  constraints on the reliability of solar and wind power in the united
  states,'' \emph{Energy \& Environmental Science}, vol.~11, no.~4, pp.
  914--925, 2018.

\bibitem{sepulveda2018role}
N.~A. Sepulveda, J.~D. Jenkins, F.~J. de~Sisternes, and R.~K. Lester, ``The
  role of firm low-carbon electricity resources in deep decarbonization of
  power generation,'' \emph{Joule}, vol.~2, no.~11, pp. 2403--2420, 2018.

\bibitem{lee2022targeted}
K.~Lee, X.~Geng, S.~Sivaranjani, B.~Xia, H.~Ming, S.~Shakkottai, and L.~Xie,
  ``Targeted demand response for mitigating price volatility and enhancing grid
  reliability in synthetic texas electricity markets,'' \emph{Iscience},
  vol.~25, no.~2, p. 103723, 2022.

\bibitem{doe}
\BIBentryALTinterwordspacing
{U.S. Department of Energy}, ``Long duration storage shot,'' 2021. [Online].
  Available: \url{https://www.energy.gov/eere/long-duration-storage-shot}
\BIBentrySTDinterwordspacing

\bibitem{zakeri2015electrical}
B.~Zakeri and S.~Syri, ``Electrical energy storage systems: A comparative life
  cycle cost analysis,'' \emph{Renewable and sustainable energy reviews},
  vol.~42, pp. 569--596, 2015.

\bibitem{safaei2015much}
H.~Safaei and D.~W. Keith, ``How much bulk energy storage is needed to
  decarbonize electricity?'' \emph{Energy \& Environmental Science}, vol.~8,
  no.~12, pp. 3409--3417, 2015.

\bibitem{jenkins2021long}
J.~D. Jenkins and N.~A. Sepulveda, ``Long-duration energy storage: A blueprint
  for research and innovation,'' \emph{Joule}, vol.~5, no.~9, pp. 2241--2246,
  2021.

\bibitem{sepulveda2021design}
N.~A. Sepulveda, J.~D. Jenkins, A.~Edington, D.~S. Mallapragada, and R.~K.
  Lester, ``The design space for long-duration energy storage in decarbonized
  power systems,'' \emph{Nature Energy}, vol.~6, no.~5, pp. 506--516, 2021.

\bibitem{kucukali2014finding}
S.~Kucukali, ``Finding the most suitable existing hydropower reservoirs for the
  development of pumped-storage schemes: An integrated approach,''
  \emph{Renewable and Sustainable Energy Reviews}, vol.~37, pp. 502--508, 2014.

\bibitem{blanco2018review}
H.~Blanco and A.~Faaij, ``A review at the role of storage in energy systems
  with a focus on power to gas and long-term storage,'' \emph{Renewable and
  Sustainable Energy Reviews}, vol.~81, pp. 1049--1086, 2018.

\bibitem{fertig2011economics}
E.~Fertig and J.~Apt, ``Economics of compressed air energy storage to integrate
  wind power: A case study in ercot,'' \emph{Energy Policy}, vol.~39, no.~5,
  pp. 2330--2342, 2011.

\bibitem{guerra2020value}
O.~J. Guerra, J.~Zhang, J.~Eichman, P.~Denholm, J.~Kurtz, and B.-M. Hodge,
  ``The value of seasonal energy storage technologies for the integration of
  wind and solar power,'' \emph{Energy \& Environmental Science}, vol.~13,
  no.~7, pp. 1909--1922, 2020.

\bibitem{conejo2005day}
A.~J. Conejo, M.~A. Plazas, R.~Espinola, and A.~B. Molina, ``Day-ahead
  electricity price forecasting using the wavelet transform and arima models,''
  \emph{IEEE transactions on power systems}, vol.~20, no.~2, pp. 1035--1042,
  2005.

\bibitem{kirschen2018fundamentals}
D.~S. Kirschen and G.~Strbac, \emph{Fundamentals of power system
  economics}.\hskip 1em plus 0.5em minus 0.4em\relax John Wiley \& Sons, 2018.

\bibitem{riddervold2021internal}
H.~O. Riddervold, E.~K. Aasg{\aa}rd, L.~Haukaas, and M.~Korp{\aa}s, ``Internal
  hydro-and wind portfolio optimisation in real-time market operations,''
  \emph{Renewable Energy}, vol. 173, pp. 675--687, 2021.

\bibitem{li2014economic}
N.~Li and K.~W. Hedman, ``Economic assessment of energy storage in systems with
  high levels of renewable resources,'' \emph{IEEE Transactions on Sustainable
  Energy}, vol.~6, no.~3, pp. 1103--1111, 2014.

\bibitem{o2017efficient}
C.~O’Dwyer, L.~Ryan, and D.~Flynn, ``Efficient large-scale energy storage
  dispatch: challenges in future high renewable systems,'' \emph{IEEE
  Transactions on Power Systems}, vol.~32, no.~5, pp. 3439--3450, 2017.

\bibitem{chen2022battery}
Y.~Chen and R.~Baldick, ``Battery storage formulation and impact on day ahead
  security constrained unit commitment,'' \emph{IEEE Transactions on Power
  Systems}, 2022.

\bibitem{huang2021configuration}
B.~Huang, Y.~Chen, and R.~Baldick, ``A configuration based pumped storage hydro
  model in the miso day-ahead market,'' \emph{IEEE Transactions on Power
  Systems}, vol.~37, no.~1, pp. 132--141, 2021.

\bibitem{nyiso_rtd}
\BIBentryALTinterwordspacing
``Nyiso manual 12 transmission and dispatch operations manual,'' 2021.
  [Online]. Available:
  \url{https://www.nyiso.com/documents/20142/2923301/trans_disp.pdf/9d91ad95-0281-2b17-5573-f054f7169551}
\BIBentrySTDinterwordspacing

\bibitem{caiso_rtd}
\BIBentryALTinterwordspacing
``California iso real-time dispatch multi-interval optimizatio,'' 2021.
  [Online]. Available:
  \url{http://www.caiso.com/Documents/EnergyStorageEnhancementsMIO-Presentation-Oct1_2021.pdf}
\BIBentrySTDinterwordspacing

\bibitem{guo2021pricing}
Y.~Guo, C.~Chen, and L.~Tong, ``Pricing multi-interval dispatch under
  uncertainty part i: Dispatch-following incentives,'' \emph{IEEE Transactions
  on Power Systems}, vol.~36, no.~5, pp. 3865--3877, 2021.

\bibitem{zhao2019multi}
J.~Zhao, T.~Zheng, and E.~Litvinov, ``A multi-period market design for markets
  with intertemporal constraints,'' \emph{IEEE Transactions on Power Systems},
  vol.~35, no.~4, pp. 3015--3025, 2019.

\bibitem{conejo2002self}
A.~J. Conejo, J.~M. Arroyo, J.~Contreras, and F.~A. Villamor, ``Self-scheduling
  of a hydro producer in a pool-based electricity market,'' \emph{IEEE
  Transactions on power systems}, vol.~17, no.~4, pp. 1265--1272, 2002.

\bibitem{ferc_order_841}
\BIBentryALTinterwordspacing
{FERC}, ``Electric storage participation in markets operated by regional
  transmission organizations and independent system operators,'' 2018.
  [Online]. Available:
  \url{https://www.ferc.gov/whats-new/comm-meet/2018/021518/E-1.pdf}
\BIBentrySTDinterwordspacing

\bibitem{sakti2018review}
A.~Sakti, A.~Botterud, and F.~O’Sullivan, ``Review of wholesale markets and
  regulations for advanced energy storage services in the united states:
  Current status and path forward,'' \emph{Energy policy}, vol. 120, pp.
  569--579, 2018.

\bibitem{caiso_es}
L.~Carr, G.~Murtaugh, J.~Powers, and B.~Sparks, ``Caiso energy storage and
  distributed energy resources phase 4 final proposal,'' 2020.

\bibitem{wang2017look}
Y.~Wang, Y.~Dvorkin, R.~Fernandez-Blanco, B.~Xu, T.~Qiu, and D.~S. Kirschen,
  ``Look-ahead bidding strategy for energy storage,'' \emph{IEEE Transactions
  on Sustainable Energy}, vol.~8, no.~3, pp. 1106--1117, 2017.

\bibitem{shafiee2016risk}
S.~Shafiee, H.~Zareipour, A.~M. Knight, N.~Amjady, and B.~Mohammadi-Ivatloo,
  ``Risk-constrained bidding and offering strategy for a merchant compressed
  air energy storage plant,'' \emph{IEEE Transactions on Power Systems},
  vol.~32, no.~2, pp. 946--957, 2016.

\bibitem{krishnamurthy2017energy}
D.~Krishnamurthy, C.~Uckun, Z.~Zhou, P.~R. Thimmapuram, and A.~Botterud,
  ``Energy storage arbitrage under day-ahead and real-time price uncertainty,''
  \emph{IEEE Transactions on Power Systems}, vol.~33, no.~1, pp. 84--93, 2017.

\bibitem{thatte2013risk}
A.~A. Thatte, L.~Xie, D.~E. Viassolo, and S.~Singh, ``Risk measure based robust
  bidding strategy for arbitrage using a wind farm and energy storage,''
  \emph{IEEE Transactions on Smart Grid}, vol.~4, no.~4, pp. 2191--2199, 2013.

\bibitem{salles2017potential}
M.~B. Salles, J.~Huang, M.~J. Aziz, and W.~W. Hogan, ``Potential arbitrage
  revenue of energy storage systems in pjm,'' \emph{Energies}, vol.~10, no.~8,
  p. 1100, 2017.

\bibitem{xu2019operational}
B.~Xu, A.~Botterud, and M.~Korpas, ``Operational valuation for energy storage
  under multi-stage price uncertainties,'' \emph{59th IEEE Conference on
  Decision and Control}, 2020.

\bibitem{castillo2013profit}
A.~Castillo and D.~F. Gayme, ``Profit maximizing storage allocation in power
  grids,'' in \emph{52nd IEEE Conference on Decision and Control}.\hskip 1em
  plus 0.5em minus 0.4em\relax IEEE, 2013, pp. 429--435.

\bibitem{nyiso}
\BIBentryALTinterwordspacing
{New York Independent System Operator, Inc. (NYISO)}, ``Energy market \&
  operational data,'' 2019. [Online]. Available:
  \url{https://www.nyiso.com/energy-market-operational-data}
\BIBentrySTDinterwordspacing

\bibitem{zheng2022arbitraing}
N.~Zheng, J.~Jaworski, and B.~Xu, ``Arbitraging variable efficiency energy
  storage using analytical stochastic dynamic programming,'' \emph{IEEE
  Transactions on Power Systems}, 2022.

\end{thebibliography}



%

\end{document}